\documentclass{svproc}

\usepackage{url}

\usepackage{amsfonts}
\usepackage{amsmath}
\usepackage{tikz}
\usepackage{graphicx, subcaption}
\usepackage[english]{babel}
\usepackage{booktabs}
\usepackage{caption}
\usepackage[per-mode=symbol, group-digits=integer]{siunitx}
\usepackage[algo2e, linesnumbered, noline, noend, ruled]{algorithm2e}
\newcommand*\diff{\mathop{}\!\mathrm{d}}
\newcommand\norm[1]{\left\lVert#1\right\rVert}

\begin{document}
\mainmatter              
\title{Computing Multiple Local Minimizers for the Topology Optimization of Bipolar Plates in Electrolysis Cells}
\titlerunning{Computing Multiple Local Minimizers for Optimization of Electrolysis Cells}  
%
\author{Leon Baeck\inst{1} \and Sebastian Blauth\inst{1}
Christian Leith\"auser\inst{1} \and Ren\'e Pinnau\inst{2} \and Kevin Sturm\inst{3}}
\authorrunning{Leon Baeck et al.} 
\institute{Fraunhofer Institute for Industrial Mathematics ITWM, Kaiserslautern, Germany,\\
\email{[leon.baeck, sebastian.blauth, christian.leithaeuser]@itwm.fraunhofer.de},
\and
RPTU Kaiserslautern, Germany,\\
\email{pinnau@mathematik.uni-kl.de}
\and
TU Vienna, Vienna, Austria,\\
\email{kevin.sturm@asc.tuwien.ac.at}
}

\maketitle

\begin{abstract}
In this paper we consider the topology optimization for a bipolar plate of a hydrogen electrolysis cell. We use the Borvall-Petersson model to describe the fluid flow and derive a criterion for a uniform flow distribution in the bipolar plate. Furthermore, we introduce a novel deflation approach to compute multiple local minimizers of topology optimization problems. The approach is based on a penalty method that discourages convergence towards previously found solutions. Finally, we demonstrate this technique on the topology optimization for bipolar plates and show that multiple distinct local solutions can be found.
\end{abstract}

\section{Introduction}

Hydrogen electrolysis cells are of major importance for sustainable energy production. Typically, hydrogen electrolysis cells use (green) electrical energy to split water into oxygen and hydrogen. One special type of electrolysis cell is the proton exchange membrane (PEM) electrolysis cell. An important role for the performance of such cells plays the anode side bipolar plate, which distributes water over the cell. To increase the cell efficiency, a uniform water distribution over the whole cell is desirable. In order to achieve such a uniform flow throughout the bipolar plate, we use techniques from topology optimization based on topological sensitivity analysis \cite{Sokolowski1999Topology}. For more details on PEM electrolysis cells, we refer the reader, e.g., to \cite{Metz2023}

Topology optimization was first introduced in the context of solid mechanics in \cite{Eschenauer1994Topology}, but has since been used in a wide variety of applications, including fluid dynamics \cite{Borrvall2003Topology}, \cite{NSa2016Topological}. It considers the optimization of a shape functional by either adding or removing material from the domain. Here, we perform topology optimization based on topological sensitivity analysis. We use a level-set approach \cite{Amstutz2006new} to solve the arising topology optimization problem. The level-set function is consequently updated by the topological derivative, which measures the sensitivity of a shape functional with respect to infinitesimal topological changes \cite{Sokolowski1999Topology}.

As topology optimization problems usually can attain multiple local solutions, deflation techniques can be used to compute local minimizers that perform better globally, see e.g. \cite{Papadopoulos2021Topology}. Deflation was already applied to find multiple solutions for nonlinear PDEs in \cite{Brow1971Deflation}, \cite{Farrell2015Deflation} or in the context of density based topology optimization \cite{Papadopoulos2021Topology}. Here, we introduce a novel approach using a deflation operator in the objective function by an additional penalty term, which allows us to find different local minimizers.

This paper is structured as follows. In Section 2, we introduce a model for the flow in the anode side bipolar plate using the Borvall-Petersson model from \cite{Borrvall2003Topology}. Additionally, we present a criterion for a uniform flow distribution that we first proposed in \cite{Baeck2023Topology}. We state the corresponding topology optimization problem and its topological derivative. In Section 3, we present a deflation approach to compute multiple local solutions of topology optimization problems. Finally, in Chapter 4, we numerically solve the topology optimization problem and present our numerical results. Particularly, we focus on using the previously introduced deflation technique to compute multiple local solutions. The obtained results show that our deflation technique can be used to create novel bipolar plate designs.

\section{Model Problem for the Anode Side Bipolar Plate}

We begin with introducing the topology optimization problem for a uniform flow distribution in the bipolar plate and state its topological derivative. For more details we refer to our previous work \cite{Baeck2023Topology}.

\subsection{The Borvall-Petersson Model}

First, we assume to have an open and bounded domain $D\subset\mathbb{R}^d$, the so-called hold all domain. We restrict ourselves to $d=2$ here, but mention that the following model has a three dimensional interpretation \cite{Borrvall2003Topology}. We identify the open and measurable set $\Omega\subset D$ as the fluid region and its complement $D\setminus\bar{\Omega}$ as the solid part of the hold all domain $D$. The fluid flow is then modeled by the Borvall Petersson model \cite{Borrvall2003Topology}. We remark that, for the sake of simplicity, we consider Stokes instead of Navier-Stokes equation here. We assume that the boundary $\Gamma=\partial D$ of the hold-all domain $D$ is divided into three parts: The inlet $\Gamma_{\mathrm{in}}$, where the water enters the plate with velocity $u_{\mathrm{in}}$, the boundary $\Gamma_{\mathrm{wall}}$, where we have a no-slip condition, and the outlet $\Gamma_{\mathrm{out}}$, where the water exits the plate with velocity $u_{\mathrm{out}}$. The non-dimensional model reads
\begin{equation}
\begin{array}{r c l l}
-\Delta u+ \alpha u + \nabla p &=& 0 &\mathrm{in }\  D,\\
\mathrm{div}(u) &=& 0 &\mathrm{in }\ D, \\
u &=& u_{\mathrm{in}}  &\mathrm{on }\ \Gamma_{\mathrm{in}}, \\
u &=& u_{\mathrm{out}} &\mathrm{on }\ \Gamma_{\mathrm{out}}, \\
u &=& 0  &\mathrm{on }\ \Gamma_{\mathrm{wall}}, \\
\int_Dp\diff x &=& 0,
\end{array}
\label{stokesdarcy}
\end{equation}
where $u$ and $p$ denote the fluid velocity and pressure, respectively. The last constraint is needed to guarantee uniqueness of the pressure $p$. We choose the outflow profile $u_{\mathrm{out}}$ such that it is compatible to the inflow profile $u_{\mathrm{out}}$. To distinguish between the solid and fluid part of the domain the inverse permeability $\alpha$ is used, which is given by
\begin{equation*}
\alpha(x) = \begin{cases} 
\alpha_U &\mathrm{if}\ x\in D\setminus\bar{\Omega},\\
\alpha_L &\mathrm{if}\ x\in\Omega,
\end{cases}
\end{equation*}
where $\alpha_L$ and $\alpha_U$ are positive constants which will be given explicitly in Chapter $4$. In particular, $\alpha$ is chosen small inside the fluid region and large in the solid part. For a detailed derivation of the model we refer to \cite{Borrvall2003Topology}.

\subsection{Uniform Flow Distribution}

To increase the cell efficiency, it is important that the water is distributed uniformly over the cell. For that we introduce a threshold velocity magnitude $u_t > 0$ that should be reached by the flow in the fluid part $\Omega$ in order to avoid dead spots. We extend this target for the fluid velocity to the hold all domain $D$ by using the smoothing characteristics of the heat equation to obtain a smoothed velocity $u_s$. We discretize the heat equation by one implicit Euler step with step length $\Delta t$ and arrive at
\begin{equation}
\begin{array}{r c l l}
\frac{u_s-u}{\Delta t}-\Delta u_s &=& 0 & \mathrm{in}\ D,\\
\nabla u_s\cdot n &=& 0 & \mathrm{on}\ \Gamma,
\end{array}
\label{heat1step}
\end{equation}
where $n$ denotes the outer unit normal vector on $\Gamma$. Using this smoothed velocity, we arrive at the criterion
\begin{equation}
\label{optgoal2}
\norm{u_s(x)}\geq u_t \ \mathrm{for\ all}\ x\in D,
\end{equation}
where $\norm{\cdot}$ denotes the Euclidean norm on $\mathbb{R}^d$. For the sake of brevity, we only gave a short explanation here. For more details we refer to our previous work \cite{Baeck2023Topology} and Remark $\ref{remark}$. Using a Moreau-Yosida regularization \cite{Hinze2009Optimization} of $(\ref{optgoal2})$ yields the following cost functional
\begin{equation}
\label{objective}
J(\Omega,u) = \int_{D} \min \left( 0, \norm{u_s} - u_t \right)^2 \diff x.
\end{equation}
\begin{remark}
	\label{remark}
	The main idea of the constraint $(\ref{optgoal2})$ is to ensure that each part of the cell receives a sufficiently large flow and that no dead spots occur. The reason of using the smoothed velocity $u_s$ for the target velocity goal is: Large solid inclusions degrade the cell efficiency as they restrict the flow in the underlying porous transport layer. Small obstacles, on the other hand, do not hinder this flow and thus are desirable. We expect that small obstacles will not be ``resolved'' by the smoothed velocity $u_s$ for a suitable time step $\Delta t$ and thus the target velocity goal $(\ref{optgoal2})$ will be reached. For large solid inclusions on the contrary, we predict that they are still ``resolved'' by the smoothed velocity meaning that the target velocity goal $(\ref{optgoal2})$ will not be reached here. Consequently, we do not expect them to appear in the final shapes. Again, we refer the reader to \cite{Baeck2023Topology} for more details.
\end{remark}

\subsection{Topology Optimization Problem}

We summarize formulas $(\ref{stokesdarcy})$, $(\ref{heat1step})$ and $(\ref{objective})$ to arrive at our topology optimization problem
\begin{equation}
\begin{array}{l l}
\min_{\Omega}J(\Omega,u)&= \int_{D} \min \left( 0, \norm{u_s} - u_t \right)^2 \diff x \\
&\mathrm{s.t.}\ (\ref{stokesdarcy})\ \mathrm{and}\ (\ref{heat1step})\ \mathrm{hold}, \\
&V_L \leq |\Omega| \leq V_U.
\end{array}
\label{optproblem}
\end{equation}
Additionally, we introduce a volume constraint with lower bound $V_L$ and upper bound $V_U$ in order to gain nontrivial results. We want to perform a topology optimization for problem $(\ref{optproblem})$ using techniques from topological sensitivity analysis, which we briefly present in the subsequent section.

\subsection{Topological Sensitivity Analysis}

We give a brief definition of the topological derivative, see e.g. \cite{Sokolowski1999Topology}. The idea of the topological derivative is to measure the sensitivity of a shape functional $\Omega\in\mathcal{A}\mapsto S(\Omega)\in\mathbb{R}$, with $\mathcal{A}=\{\Omega\subset D\ \mathrm{s.t.}\ \Omega \ \mathrm{open}\}$, with respect to infinitesimal topological changes. We introduce the perturbation $\Omega_{z,\epsilon}$ of $\Omega$ around $z\in D\setminus\partial\Omega$ by
\begin{equation*}
\Omega_{z,\epsilon}=\begin{cases}
\Omega\setminus\bar{\omega}_{z,\epsilon}, & z\in\Omega\\
\Omega\cup\omega_{z,\epsilon}, & z\in D\setminus\bar{\Omega},
\end{cases}
\end{equation*}
where $\omega_{z,\epsilon}=z+\epsilon\omega$ with $\omega\subset\mathbb{R}^d$ being a simply connected domain with $0\in\omega$.  Furthermore, for a positive function $l$ with $\lim_{\epsilon\searrow0}l(\epsilon)=0$, the topological derivative is defined as
\begin{equation}
D_TS(\Omega,\omega)(z)=\lim_{\epsilon\searrow0}\frac{S(\Omega_{z,\epsilon})-S(\Omega)}{l(\epsilon)},
\end{equation}
if the limit exists. For numerical purposes we also state the generalized topological derivate
\begin{equation*}
	\mathcal{D}_TS(\Omega,\omega)(z)=\begin{cases}
	-D_TS(\Omega,\omega)(z), & z\in\Omega\\
	+D_TS(\Omega,\omega)(z), & z\in D\setminus\bar{\Omega}.
	\end{cases}
\end{equation*}

The generalized topological derivative of our model problem $(\ref{optproblem})$ reads
\begin{equation}
\label{topder}
\mathcal{D}_TJ(z)=-(\alpha_U-\alpha_L)u(z)v(z)
\end{equation}
for all $z\in D\setminus\partial\Omega$. Here, $u$ is the weak solution of $(\ref{stokesdarcy})$ and $v$ is the adjoint velocity which solves
\begin{equation}
\begin{array}{r c l l}
-\Delta v+ \alpha v + \nabla q-\frac{1}{\Delta t}v_s &=& 0 & \mathrm{in}\ D,\\
\mathrm{div}(v) &=& 0 & \mathrm{in}\ D, \\
v &=& 0 & \mathrm{on}\ \Gamma, \\
\int_Dq\diff x &=& 0,
\end{array}
\label{stokesadj}
\end{equation}
and $v_s$ is the adjoint of the smoothed velocity which solves the equation
\begin{equation}
\label{heatadjoint}
\begin{array}{r c l l}
\frac{1}{\Delta t}v_s-\Delta v_s &=& 2\frac{u_s}{\norm{u_s}}\min(0,\norm{u_s}-u_t) & \mathrm{in}\ D,\\
\nabla v_s\cdot n &=& 0 & \mathrm{on}\ \partial D.
\end{array}
\end{equation}
A rigorous computation of the topological derivative, e.g. with an averaged adjoint approach, see e.g. \cite{Sturm2020Topology}, is beyond the scope of this paper and a topic of future research.

\section{Deflation}

Typically, topology optimization problems can attain multiple local solutions. To gain local minimizers that perform better globally, deflation techniques have already been used in the context of density based topology optimization, see e.g. \cite{Papadopoulos2021Topology}. Here, we introduce a novel deflation technique for gradient-based topology optimization. Assuming we already found a local minimizer, we alter the objective function of a topology optimization problem by adding a penalty term which penalizes shapes that are close to previously found local solutions. Additionally, the penalty term should vanish if we are far from the previously found local minimizers. To measure the distance between two shapes $\Omega_1$ and $\Omega_2$, we consider the distance between the two corresponding characteristic functions $\chi_1$ and $\chi_2$ in the $L^2$-sense, i.e.
\begin{equation}
	\mathrm{dist}(\Omega_1,\Omega_2)=\norm{\chi_1-\chi_2}_{L^2(D)}.
	\label{distance}
\end{equation}
Then, we propose the following penalty function
\begin{equation}
P_{\gamma,\delta}(\Omega_1,\Omega_2)=\mathrm{exp}\left(\delta\left(\frac{\gamma^2}{\norm{\chi_1-\chi_2}_{L^2(D)}}-\norm{\chi_1-\chi_2}_{L^2(D)}\right)\right),
\label{penalty}
\end{equation}
where $\gamma$ and $\delta$ are positive constants. The parameter $\gamma$ controls the similarity of the two shapes. It is straightforward to see that this penalty function indeed increases when the two shapes get close to each other
\begin{equation*}
	\lim_{\mathrm{dist}(\Omega_1,\Omega_2)\rightarrow0}P_{\gamma,\delta}(\Omega_1,\Omega_2)=+\infty.
\end{equation*}
On the other hand, we know that $\norm{\chi_1-\chi_2}_{L^2(D)}$ is bounded from above by the Lebesgue measure of the hold all domain $D$. Thus, we can not expect a vanishing penalty function in case that the two shapes $\Omega_1$ and $\Omega_2$ are far from each other. To still guarantee that the penalty function gets close to zero, we choose a large value for the constant $\delta$.

This procedure can then be applied iteratively to compute even more local solutions. Assuming we already found local minimizers $\Omega_0,...\Omega_i$ with characteristic functions $\chi_0,...,\chi_i$, we consider the penalty function
\begin{equation}
P_{i,\gamma,\delta}\left(\Omega,\Omega_0,...,\Omega_i\right)=\sum_{j=0}^{i}\mathrm{exp}\left(\delta\left(\frac{\gamma^2}{\norm{\chi-\chi_j}_{L^2(D)}}-\norm{\chi-\chi_j}_{L^2(D)}\right)\right).
\label{penaltylarge}
\end{equation}

Turning back to our topology optimization problem $(\ref{optproblem})$, we add the previous penalty term $(\ref{penaltylarge})$ to the objective and arrive at
\begin{equation}
\begin{array}{l l}
\min_{\Omega}\tilde{J}_i(\Omega,u)&= \int_{D} \min \left( 0, \norm{u_s} - u_t \right)^2 \diff x + P_{i,\gamma,\delta}\left(\Omega,\Omega_0,...,\Omega_i\right) \\
&\mathrm{s.t.}\ (\ref{stokesdarcy})\ \mathrm{and}\ (\ref{heat1step})\ \mathrm{hold}, \\
&V_L \leq |\Omega| \leq V_U.
\end{array}
\label{optproblemdeflation}
\end{equation}
In order to perform topology optimization for $(\ref{optproblemdeflation})$, we state the corresponding generalized topological derivative
\begin{equation}
\begin{aligned}
\mathcal{D}_T\tilde{J}_i(\Omega)&=\mathcal{D}_TJ(\Omega)+\mathcal{D}_TP_{i,\gamma,\delta}\left(\Omega,\Omega_0,...,\Omega_i\right) \\
&=-(\alpha_U-\alpha_L)u(z)v(z) \\
&\phantom{=}-\sum_{j=0}^{i}\delta\left(\frac{\gamma^2}{2\norm{\chi-\chi_j}^3}+\frac{1}{2\norm{\chi-\chi_j}}\right)(1-2\chi)\\
&\phantom{=-}\mathrm{exp}\left(\delta\left(\frac{\gamma^2}{\norm{\chi-\chi_j}_{L^2(D)}}-\norm{\chi-\chi_j}_{L^2(D)}\right)\right).
\end{aligned}
\end{equation}
As before, $u$ is the weak solution $(\ref{stokesdarcy})$ and v solves $(\ref{stokesadj})$ in the weak sense. The computation of the topological derivative of the penalty function is straightforward. Applying the definition of the topological derivative to the distance of two shapes $(\ref{distance})$ and straightforward calculations yield the above formula. The procedure for computing multiple local solutions of $(\ref{optproblem})$ is described in Algorithm $\ref{algo:deflation}$. In Step $4$ we perform a restart with problem $(\ref{optproblem})$ with the previously computed deflated solution as the starting shape. This is done to guarantee that we reach a solution of the actual topology optimization problem $(\ref{optproblem})$ and to omit the disturbance of the penalty function. Hence, our deflation approach can be interpreted as a technique to generate distinct initial guesses for the optimization in a systematic way.

\begin{algorithm2e}[!t]
	\KwIn{Parameters $\delta>0$ and $\gamma>0$ for the penalty function, number of desired local minimizers $n$}
	Compute a minimizer $\Omega_0$ of $(\ref{optproblem})$ and store its characteristic function $\chi_0$\\
	\For{$i=0,1,2,\dots,n-1$}{
		Compute a minimizer $\tilde{\Omega}_{i+1}$ of $(\ref{optproblemdeflation})$ \\
		Compute a solution $\Omega_{i+1}$ of $(\ref{optproblem})$ with starting shape $\tilde{\Omega}_{i+1}$ \\
		Store the characteristic function $\chi_{i+1}$ of $\Omega_{i+1}$ \\
	}
	\caption{Deflation procedure for computing multiple local minimizers of $(\ref{optproblem})$.}
	\label{algo:deflation}
\end{algorithm2e}

\section{Numerical Results}

\begin{figure}
	\centering
	\begin{tikzpicture}[scale=0.3]
	\draw[line width = 0.35mm] (0,0) -- (10,0) -- (10,10) -- (0,10) -- (0,0);
	\draw[black, line width = 0.30mm]   plot[smooth,domain=3.5:6.5] ({-(\x-3.5)*(6.5-\x)}, \x);
	\draw[black, line width = 0.30mm]   plot[smooth,domain=3.5:6.5] ({(\x-3.5)*(6.5-\x)+10}, \x);
	\draw[black, line width = 0.30mm, ->] (-1.6875,5.75) -- (0,5.75);
	\draw[black, line width = 0.30mm, ->] (-2.25,5) -- (0,5);
	\draw[black, line width = 0.30mm, ->] (-1.6875,4.25) -- (0,4.25);
	\draw[black, line width = 0.30mm, ->] (10,5.75) -- (11.6875,5.75);
	\draw[black, line width = 0.30mm, ->] (10,5) -- (12.25,5);
	\draw[black, line width = 0.30mm, ->] (10,4.25) -- (11.6875,4.25);
	\draw[black, line width = 0.30mm, <->] (0.5,5) -- (0.5,10);
	\draw[black, line width = 0.30mm, <->] (0.5,0) -- (0.5,3.5);
	\draw[black, line width = 0.30mm, <->] (0,10.5) -- (10,10.5);
	\draw[black, line width = 0.30mm, <->] (13,0) -- (13,10);
	\node[] (a) at (5,5) {$D$};
	\node[] (b) at (15.25,5) {$l=1.0$};
	\node[] (b) at (5,11.25) {$l=1.0$};
	\node[] (b) at (1.75,1.75) {$0.35$};
	\node[] (b) at (1.5,7.5) {$0.5$};
	\end{tikzpicture}
	\caption{Schematic setup of the hold all domain $D$.}
	\label{fig_model}
\end{figure}
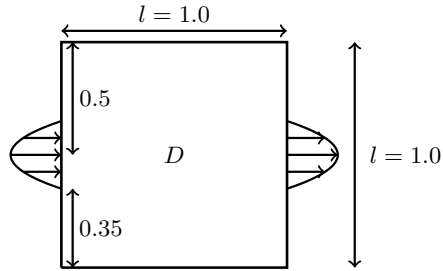

The goal is to use the deflation approach to compute local minimizers of $(\ref{optproblem})$ with algorithm \ref{algo:deflation}. We start by introducing the setting for the numerical experiments. The setup of the setting is displayed in Figure $\ref{fig_model}$, where the hold all domain $D=(0,1)\times(0,1)$ is the two dimensional unit square. The inflow profile $u_\mathrm{in}$, which is applied on the inflow boundary $\Gamma_{\mathrm{in}}$, reads
\begin{equation*}
u_{\mathrm{in}}=\begin{bmatrix}
-\frac{400}{9}(y-\frac{7}{20})(y-\frac{13}{20})\\
0
\end{bmatrix} \; \mathrm{for}\ x=0 \ \mathrm{and} \ \frac{7}{20}\leq y\leq \frac{13}{20}.
\end{equation*}
Analogously, we introduce the outflow profile $u_\mathrm{out}$ on the outflow boundary $\Gamma_{\mathrm{out}}$ by
\begin{equation*}
u_{\mathrm{out}}=\begin{bmatrix}
-\frac{400}{9}(y-\frac{7}{20})(y-\frac{13}{20})\\
0
\end{bmatrix} \; \mathrm{for}\ x=1 \ \mathrm{and} \ \frac{7}{20}\leq y\leq \frac{13}{20}.
\end{equation*}
For the inverse permeability we choose the values 
\begin{equation*}
\alpha_L=\frac{2.5}{100^2}, \;\; \alpha_U=\frac{2.5}{0.0025^2}
\end{equation*}
to avoid fluid flow through solid, which was observed in numerical tests for smaller values of $\alpha_U$, e.g. for the choices in \cite{Borrvall2003Topology}. The target velocity is chosen as $u_t=0.1$ and the lower and upper values for the volume constraint are given by $V_L=0.5$ and $V_U=0.7$, respectively. Finally, we set the step length for the smoothing equation to $\Delta t=10^{-3}$ and we refer to \cite{Baeck2023Topology} for a discussion of the step length.

\subsection{Implementation}

To perform the topology optimization, we use a gradient-based solution algorithm which was introduced in \cite{Amstutz2006new}. The main idea of the algorithm is to characterize the solid and fluid parts of the domain $D$ by a level-set function and then consequently update this level-set function using the generalized topological derivative \cite{Amstutz2006new}. The algorithm is stopped when the angle between the level-set function and the generalized topological derivative becomes smaller than a certain numerical tolerance. We set this tolerance to $\epsilon_\theta=0.035$, which is equivalent to an angle of $2$ degrees. For more details the reader is referred to \cite{Amstutz2006new}. Additionally, we refer to \cite{Blauth2023Quasi} for an overview over established topology optimization algorithms as well as novel quasi-Newton methods for topology optimization.

We give a short overview over some implementational aspects. We divide our domain in 10.000 uniform triangular elements. All underlying PDEs are solved with the software package FEniCS \cite{Alnes2015FEniCS}. For the Stokes-Darcy equation $(\ref{stokesdarcy})$ we use Taylor-Hood finite elements, and for the smoothing equation $(\ref{heat1step})$ continuous quadratic Lagrange elements.

The adjoint equations are computed with the software package cashocs \cite{Blauth2023Cashocs}, which implements a discretization of the continuous adjoint approach and thus allows for automatic adjoint computation. The software package cashocs is open source and based on FEniCS. It can be used to solve shape optimization, optimal control and topology optimization problems in an automated fashion.

The volume constraint is handled by a projection of the level-set function onto the admissible set, for more details we refer to \cite{Baeck2023Topology}.

\subsection{First Result}

A local minimizer of $(\ref{optproblem})$ is displayed in the first row of Table \ref{tab:result}, where the optimized geometry, the corresponding flow velocity norm and the smoothed flow velocity norm are shown. The velocity plots are given in a logarithmic scale. The algorithm took 53 iterations to reach the numerical target accuracy $\epsilon_\theta$. The target velocity constraint $(\ref{optgoal2})$ is fulfilled on $76.4$ percent of the hold all domain $D$. To do so the optimizer created a structure which consists of four channels.

\begin{table}[!t]
	\newcommand{\sizebox}{0.2775\textwidth}
	\newcolumntype{C}{>{\centering\arraybackslash}m{\sizebox}}
	\setlength{\tabcolsep}{0pt}
	\caption{Different local minimizers of the topology optimization problem $(\ref{optproblem})$.}
	\label{tab:result}
	\begin{tabular}{m{0.15\textwidth} @{\hskip 0.5em} C C C }
		\toprule
		Fulfillment of $(\ref{optgoal2})$ for $u_s$ & Final shape & Norm of the velocity field & Norm of the smoothed velocity field \\ 
		\midrule
		$76.4\%$ & \includegraphics[width=\sizebox]{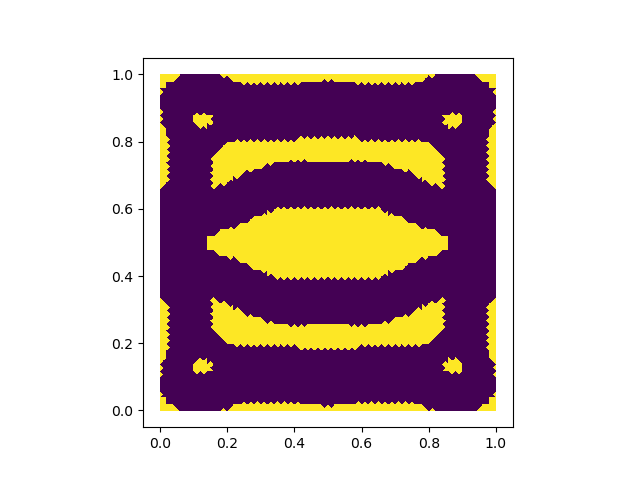} & \includegraphics[width=\sizebox]{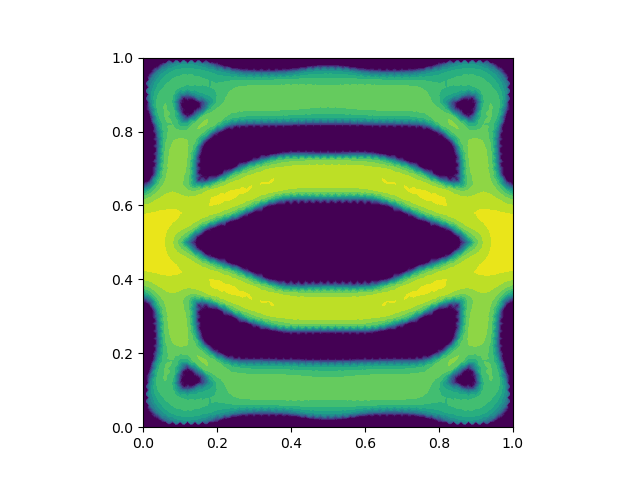} & \includegraphics[width=\sizebox]{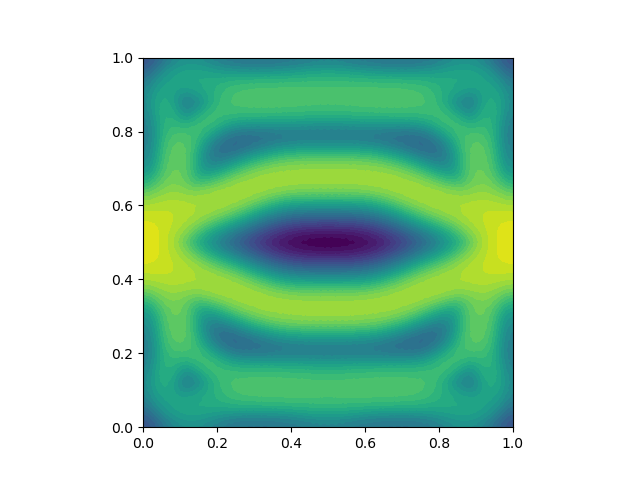} \\ 
		$90.68\%$ & \includegraphics[width=\sizebox]{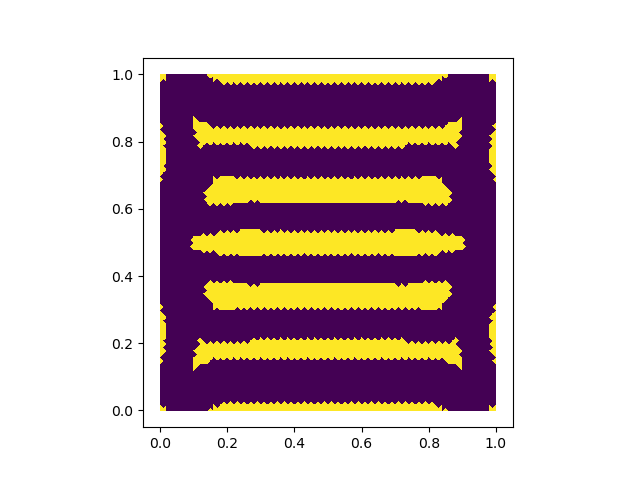} & \includegraphics[width=\sizebox]{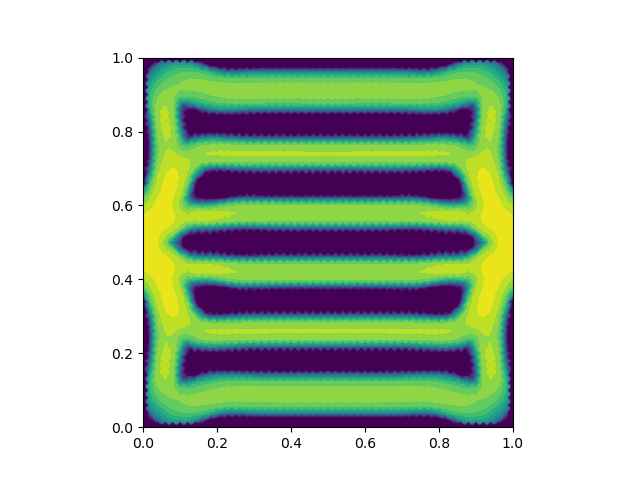} & \includegraphics[width=\sizebox]{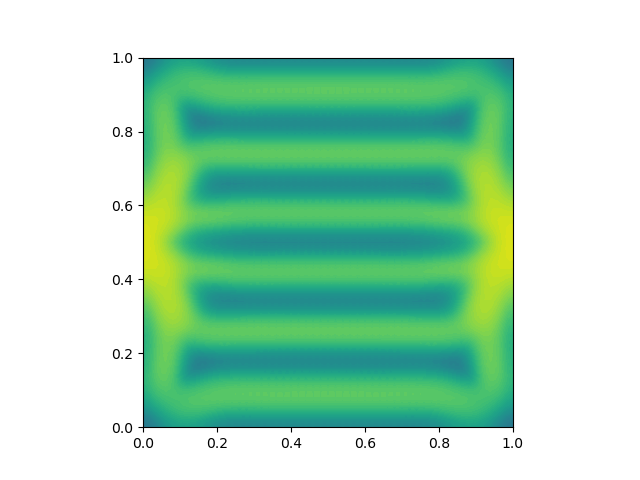} \\ 
		$98.88\%$ & \includegraphics[width=\sizebox]{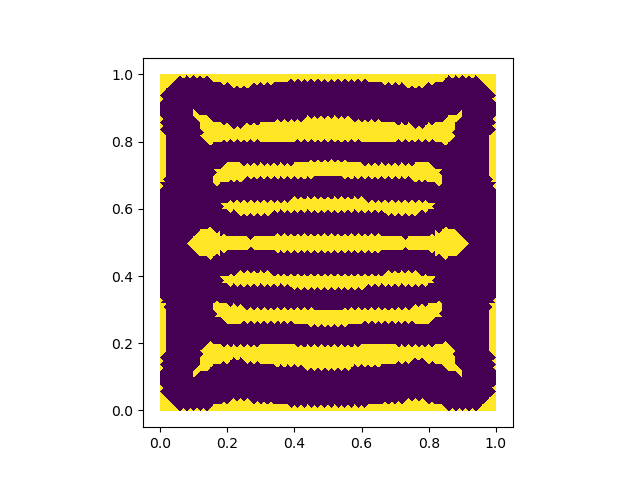} & \includegraphics[width=\sizebox]{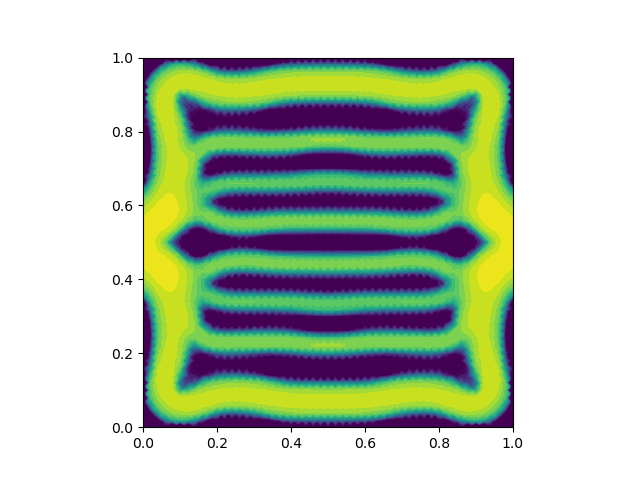} & \includegraphics[width=\sizebox]{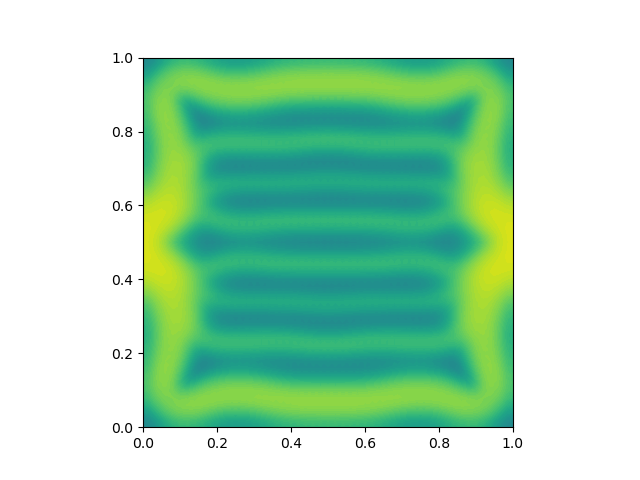} \\ 
		\bottomrule 
	\end{tabular}
\end{table}

\subsection{Deflation}

As already mentioned before, topology optimization problems tend to posess multiple local solutions. Thus, we apply our deflation technique to compute other local minimizers of the optimization problem $\ref{optproblem}$. We solve problem $(\ref{optproblemdeflation})$ multiple times, in the hope of finding local minimizers that perform better globally. For the penalty function we choose $\gamma=0.4$ and $\delta=50$. As described previously in Algorithm \ref{algo:deflation}, after solving the deflated problem $(\ref{optproblemdeflation})$, we solve the original problem $(\ref{optproblem})$ with the solution of the deflated problem as starting shape.

The first additional local solution is displayed in the second row of Table \ref{tab:result}. The optimization algorithm performed 46 iterations for the deflated and 51 iterations for the unperturbed problem. The shape performs better in the context of the target optimization goal $(\ref{optgoal2})$, as it is fulfilled on $90.68\%$ of the domain $D$. Additionally, it is apparent that the reached fluid distribution between the six channels is more uniform compared to the first result in the first row of Table \ref{tab:result}.

Repeating the described procedure, we compute a third local minimizer of problem $(\ref{optproblem})$, where 38 iterations for the deflated problem and 60 iterations for the undeflated problem were needed. The solution can be found in the third row of Table \ref{tab:result}. Again, the deflation approach produces a local solution that performs better globally, with a $98.88\%$ fulfillment of the velocity constraint $(\ref{optgoal2})$. For that the optimizer formed a structure consisting of 8 channels with a uniform flow distribution between them.

\section{Conclusion}

Summarizing, we presented a model for achieving a uniform flow distribution in the anode side bipolar plate of a hydrogen electrolysis cell. Additionally, we introduced a deflation technique for derivative-based topology optimization using a penalty approach to find multiple local solutions. We applied this approach numerically and computed several minimizers for the presented topology optimization problem. As desired, the deflation approach indeed delivered more optimal solutions. Furthermore, the numerical experiments showed that a uniform flow distribution is indeed achieved in the optimized geometries.

\end{document}